\documentclass{amsart}
\usepackage{amsmath, amsthm, amsfonts, amssymb, txfonts, graphicx, hyperref}
\newtheorem{thrm}{Theorem}[section]

\newtheorem*{cnjctr}{Conjecture}

\numberwithin{sblmm}{thrm}
\numberwithin{equation}{section}
\newtheoremstyle{named}{}{}{\itshape}{}{\bfseries}{.}{.5em}{\thmnote{#3}}
\theoremstyle{named}
\newtheorem*{NamedTheorem}{Theorem}

\begin{document}
\title{Bounded length intervals containing two primes and an almost-prime}
\author{James Maynard}
\address{Mathematical Institute, 24–-29 St Giles', Oxford, OX1 3LB}
\email{maynard@math.ox.ac.uk}
\thanks{Supported by EPSRC Doctoral Training Grant EP/P505216/1 }
\date{}
\begin{abstract}
Goldston, Pintz and Y\i ld\i r\i m have shown that if the primes have `level of distribution' $\theta$ for some $\theta>1/2$ then there exists a constant $C(\theta)$, such that there are infinitely many integers $n$ for which the interval $[n,n+C(\theta)]$ contains two primes. We show under the same assumption that for any integer $k\ge 1$ there exists constants $D(\theta,k)$ and $r(\theta,k)$, such that there are infinitely many integers $n$ for which the interval $[n,n+D(\theta,k)]$ contains two primes and $k$ almost-primes, with all of the almost-primes having at most $r(\theta,k)$ prime factors. If $\theta$ can be taken as large as $0.99$, and provided that numbers with $2$, $3$, or $4$ prime factors also have level of distribution $0.99$, we show that there are infinitely many integers $n$ such that the interval $[n,n+90]$ contains 2 primes and an almost-prime with at most 4 prime factors.
\end{abstract}
\maketitle
\section{Introduction}
We are interested in trying to understand how small gaps between primes can be. If we let $p_n$ denote the $n^{\text{th}}$ prime, it is conjectured that
\begin{equation}
\liminf_n p_{n+1}-p_n=2.
\end{equation}
This is the famous twin prime conjecture. Unfortunately we appear unable to prove any results of this strength. The best unconditional result is due to Goldston, Pintz and Y\i ld\i r\i m \cite{GPY} which states that
\begin{equation}
\liminf_n \frac{p_{n+1}-p_n}{\sqrt{\log{p_n}}(\log\log{p_n})^2}<\infty.
\end{equation}
Therefore we do not know that $\liminf p_{n+1}-p_n$ is finite.

The method of \cite{GPY} relies heavily on results about primes in arithmetic progressions. We say that the primes have `level of distribution' $\theta$ if for any constant $A$ there is a constant $C=C(A)$ such that
\begin{equation}
\sum_{q\le x^\theta(\log{x})^{-C}}\max_{\substack{a\\(a,q)=1}}\left|\sum_{\substack{p\equiv a \pmod{q}\\ p\le x}}1-\frac{\text{Li}(x)}{\phi(q)}\right|\ll_A\frac{x}{(\log{x})^A}.
\end{equation}
The Bombieri-Vinogradov theorem states that the primes have level of distribution $1/2$, and this is a major ingredient in the proof of the Goldston-Pintz-Y\i ld\i r\i m result.

If we could improve the Bombieri-Vinogradov theorem to show that the primes have level of distribution $\theta$ for some constant $\theta>1/2$, then it would follow from \cite{GPYI}[Theorem 1] that there is a constant $D=D(\theta)$ such that
\begin{equation}\liminf_n p_{n+1}-p_n<D,\end{equation}
and so there would be infinitely many bounded gaps between primes. It is believed that such improvements to the Bombieri-Vinogradov theorem are true, and Elliott and Halberstam \cite{ElliottHalberstam} conjectured the following much stronger result.
\begin{cnjctr}[Elliott-Halberstam Conjecture]
For any fixed $\epsilon>0$, the primes have level of distribution $1-\epsilon$.
\end{cnjctr}
Friedlander and Granville \cite{FriedlanderGranville} have shown that the primes do not have level of distribution $1$, and so the Elliott-Halberstam conjecture represents the strongest possible result of this type.

Under the Elliott-Halberstam conjecture the Goldston-Pintz-Y\i ld\i r\i m method gives \cite{GPYI} that
\begin{equation}
\liminf_n p_{n+1}-p_n\le 16.
\end{equation}
If we consider the length of 3 or more consecutive primes, however, we are unable to prove as strong results, even under the full strength of the Elliott-Halberstam conjecture. In particular we are unable to prove that there are infinitely many intervals of bounded length that contain at least 3 primes. The Goldston-Pintz-Y\i ld\i r\i m methods can still be used, but even with the Elliot-Halberstam conjecture we are only able to prove that
\begin{equation}
\liminf_n \frac{p_{n+2}-p_n}{\log{p_n}}=0.
\end{equation}
This should be contrasted with the following conjecture.
\begin{cnjctr}[Prime $k$-tuples conjecture]
Let $\mathcal{L}=\{L_1,\dots,L_k\}$ be a set of integer linear functions whose product has no fixed prime divisor. Then there are infinitely many $n$ for which all of $L_1(n), L_2(n), \dots,L_k(n)$ are simultaneously prime.
\end{cnjctr}
By `no fixed prime divisor' above we mean that for every prime $p$ there is an integer $n_p$ such that $L_i(n_p)$ is coprime to $p$ for all $1\le i\le k$. We call such a set of linear functions \textit{admissible}.

We note that $\{n,n+2,n+6\}$ is an admissible set of linear functions, and so the prime $k$-tuples conjecture predicts that $\liminf_n p_{n+2}-p_n\le 6$ (it is easy to verify that one cannot have $p_{n+2}-p_n<6$ for $n>2$). More generally, for any constant $k>0$ the conjecture predicts that $\liminf p_{n+k}-p_n<\infty$, and so there are infinitely many intervals of bounded size containing at least $k$ primes.

At the moment the prime $k$-tuples conjecture appears beyond the techniques currently available to us. As an approximation to the conjecture, it is common to look for \textit{almost-prime} numbers instead of primes, where almost-prime indicates that the number has only a `few' prime factors.

Graham, Goldston, Pintz and Y\i ld\i r\i m \cite{GGPY} have shown that given an integer $k$, there are infinitely many intervals of bounded length (depending on $k$) containing at least $k$ integers each with exactly two prime factors. It is a classical result of Halberstam and Richert \cite{HalberstamRichert} that there are infinitely many intervals of bounded length (depending on $k$) which contain a prime and at least $k$ numbers each with at most $r$ prime factors for $r$ sufficiently large (depending on $k$).

We investigate, under the assumption that the primes have level of distribution $\theta>1/2$, whether there are infinitely many intervals of bounded length (depending on $k$) containing $2$ primes and $k$ numbers each with at most $r$ prime factors.
\section{Initial Hypotheses}
We will work with an assumption either on the distribution of primes in arithmetic progressions of level $\theta$, or a stronger assumption on numbers with exactly $r$ prime factors each of which is of a given size.

Given constants $0\le \eta_i\le\delta_i\le 1$ for $1\le i\le r$ we define
\begin{equation}
\beta_{r,\eta,\delta}(n)=\begin{cases}
1,\qquad&\text{$n=p_1p_2\dots p_r$ with $n^{\eta_i}\le p_i\le n^{\delta_i}$ for $1\le i\le r$,}\\
0,&\text{otherwise.}
\end{cases}
\end{equation}
We put
\begin{align}
\Delta(x;q,a)&=\sum_{\substack{x<p\le 2x\\p\equiv a\pmod{q}}}1-\frac{1}{\phi(q)}\sum_{x<p\le 2x}1,\\
\Delta_{r,\eta,\delta}(x;q,a)&=\sum_{\substack{x<p\le 2x\\p\equiv a\pmod{q}}}\beta_{r,\eta,\delta}(n)-\frac{1}{\phi(q)}\sum_{x<p\le 2x}\beta_{r,\eta,\delta}(n),\\
\Delta^*(x;q)&=\max_{y\le x}\max_{\substack{a\\(a,q)=1}}\left|\Delta(y;q,a)\right|,\\
\Delta_{r,\eta,\delta}^*(x;q)&=\max_{y\le x}\max_{\substack{a\\(a,q)=1}}\left|\Delta_{r,\eta,\delta}(y;q,a)\right|.
\end{align}
We can now state the two hypotheses that we will consider, the Bombieri-Vinogradov hypothesis of level $\theta$, BV($\theta$), and the generalised Bombieri-Vinogradov hypothesis of level $\theta$ for $E_r$ numbers, GBV($\theta$,r).
\begin{NamedTheorem}[Hypothesis BV($\theta$)]For every constant $A>0$ and integer $h>0$ there is a constant $C=C(A,h)$ such that if $Q\le x^\theta(\log{x})^{-C}$ then we have
\[\sum_{q\le Q}\mu^2(q)h^{\omega(q)}\Delta^*(x;q)\ll_A x(\log{x})^{-A}.\]
\end{NamedTheorem}
\begin{NamedTheorem}[Hypothesis GBV($\theta$,r)]For every constant $A>0$ and integer $h>0$ there is a constant $C=C(A,h)$ such that if $Q\le x^\theta(\log{x})^{-C}$ then uniformly for $0\le \eta_i\le \delta_i\le 1$ ($1\le i\le r$) we have
\[\sum_{q\le Q}\mu^2(q)h^{\omega(q)}\Delta^*_{r,\eta,\delta}(x;q)\ll_A x(\log{x})^{-A}.\]
\end{NamedTheorem}
We note that by standard arguments in sieve methods (see, for example, \cite{HalberstamRichert}[Lemma 3.5]) Hypothesis BV($\theta$) follows from the primes having level of distribution $\theta$.
\section{Statement of Results}
\begin{thrm}\label{thrm:MainTheorem}
Let $k\ge 1$ be an integer. Let $1/2<\theta<0.99$. Assume Hypothesis BV($\theta$) holds. Let
\[r=\frac{240k^2}{(2\theta-1)^3}.\]
Then there are infinitely many integers $n$ such that the interval $[n,n+C(k,\theta)]$ contains two primes and $k$ integers, each with at most $r$ prime factors.
\end{thrm}
\begin{thrm}\label{thrm:ElliottHalberstam}
Let $\theta\ge 0.99$, and assume Hypothesis GBV($\theta$,r) holds for $1\le r\le 4$. Then there exist infinitely many integers $n$ such that the interval $[n,n+90]$ contains two primes and one other integer with at most 4 prime factors.\end{thrm}
\section{Proof of Theorem \ref{thrm:MainTheorem}}
We consider two finite disjoint sets of integer linear functions $\mathcal{L}^{(1)}_1=\{L^{(1)}_1,\dots,L^{(1)}_{k_1}\}$ and $\mathcal{L}^{(1)}_2=\{L^{(1)}_{k_1+1},\dots,L^{(1)}_{k_1+k_2}\}$, whose union $\mathcal{L}^{(1)}=\mathcal{L}_1^{(1)}\cup\mathcal{L}^{(1)}_2$ is admissible. (We recall that a such set is admissible if for every prime $p$ there is an integer $n_p$ such that every function evaluated at $n_p$ is coprime to $p$).

We wish to show that there are infinitely many $n$ for which two of the functions from $\mathcal{L}_1^{(1)}$ take prime values at $n$, and at all of the functions from $\mathcal{L}^{(2)}$ take almost-prime values at $n$.

Since we are only interested in showing there are infinitely many such $n$, we adopt a normalisation of our linear functions, as done originally by Heath-Brown \cite{HeathBrown} which simplifies our argument. By considering $L_i(n)=L_i^{(1)}(An+B)$ for suitable constants $A$ and $B$ we may assume that the functions $L_i$ satisfy the following conditions.
\begin{enumerate}
\item The functions $L_i(n)=a_i n+b_i$ ($1\le i \le k_1+k_2$) are distinct with $a_i>0$. 
\item Each of the coefficients $a_i$ is composed of the same primes none of which divides the $b_j$. 
\item If $i\ne j$, then any prime factor of $a_i b_j-a_j b_i$ divides each of the $a_l$.
\end{enumerate}
We let $\mathcal{L}_1=\{L_1,\dots,L_{k_1}\}$ and $\mathcal{L}_2=\{L_{k_1+1},\dots,L_{k_1+k_2}\}$.

We now consider the sum
\begin{equation}
S(N;\mathcal{L}_1,\mathcal{L}_2)=\sum_{N\le n \le 2N}\left(\sum_{L\in\mathcal{L}_1}\chi_{1}(L(n))+\sum_{L\in\mathcal{L}_2}\chi_{r}(L(n))-k_2-1\right)\left(\sum_{\substack{d|\Pi(n)\\ d\le R}}\lambda_d\right)^2,
\end{equation}
where
\begin{align}
\chi_r(n)&=\begin{cases}
1,\qquad &\text{$n$ has at most $r$ prime factors}\\
0,&\text{otherwise,}\end{cases}\\
R&=N^{\theta/2}(\log{N})^{-C},\\
\Pi(n)&=\prod_{L\in\mathcal{L}_1\cup\mathcal{L}_2}L(n),
\end{align}
and the $\lambda_d$ are real numbers which we declare later. $C>0$ is a constant chosen sufficiently large so we can use the estimates of hypotheses $BV(\theta)$ or $GBV(\theta,r)$.

If we can show that $S>0$ then we know there must be at least one $n\in[N,2N]$ for which the terms in parentheses give a positive contribution to $S$. The second term in our expression for $S$ is a square, and so is always non-negative. We see that the first term in parentheses is positive only when there are at least two primes and $k_2$ numbers each with at most $r$ prime factors amongst the $L_i(n)$ ($1\le i\le k_1+k_2$). If we choose all our original functions to be of the form $L^{(1)}_i(n)=n+h_i$ (with $h_i\ge 0$) then all these integers then lie in an interval $[m,m+H]$, where $H=\max_{i} h_i$.

Thus it is sufficient to show that $S>0$ for any large $N$ to prove Theorem \ref{thrm:MainTheorem}. We can get such a bound by following a method similar to Goldston, Pintz and Y\i ld\i r\i m \cite{GPY}, which we refer to as the GPY method.

To simplify notation we put
\begin{equation}
\Lambda^2(n)=\left(\sum_{\substack{d|\Pi(n)\\d\le R}}\lambda_d\right)^2.
\end{equation}
To avoid confusion we mention that $\Lambda^2(n)$ is unrelated to the Von-Mangold function.

We expect to be able to show that $S>0$ for suitably large $k_1$ and $r$ (depending on $k_2$) when the primes have level of distribution $\theta>1/2$. This is because the original GPY method shows that for sufficiently large size of $k_1$ (depending on $k_2$ and $\epsilon$) we can choose the $\lambda_d$ to give
\begin{equation}
\sum_{N\le n\le 2N}\sum_{L\in \mathcal{L}_1}\chi_1(L(n))\Lambda^2(n)\ge (2\theta-\epsilon)\sum_{N\le n\le 2N}\Lambda^2(n).
\end{equation}
Moreover, since $\Lambda^2(n)$ is small when $\Pi(n)$ has many prime factors, we expect for sufficiently large $r$ (depending on $k_2$ and $\epsilon$) that
\begin{equation}
\sum_{N\le n\le 2N}\sum_{L\in\mathcal{L}_2}(1-\chi_r(L(n)))\Lambda^2(n)\le \epsilon\sum_{N\le n\le 2N}\Lambda^2(n).
\end{equation}
And so provided that $\theta>1/2+\epsilon$ we expect that
\begin{equation}
S\gg_\epsilon \sum_{N\le n\le 2N}\Lambda^2(n)>0.
\end{equation}
Although the method of Graham, Goldston, Pintz and Y\i ld\i r\i m allows one to estimate similar sums involving numbers with a fixed number of prime factors, these results rely on level-of-distribution results for such numbers, which we are not assuming in Theorem \ref{thrm:MainTheorem}. Instead we proceed by noting that any integer which is at most $2N$ and has more than $r$ prime factors must have a prime factor of size at most $(2N)^{1/(r+1)}$. Thus for $n\le 2N$
\begin{equation}
\chi_r(n)\ge 1-\sum_{\substack{p|n\\p\le (2N)^{1/(r+1)}}}1.
\end{equation}
Substituting this into our expression for $S$ we have
\begin{align}
S&\ge \sum_{N\le n\le 2N}\left(\sum_{L\in \mathcal{L}_1}\chi_1(L(n))-1-\sum_{L\in \mathcal{L}_2}\sum_{\substack{p|L(n)\\p\le (2N)^{1/r+1}}}1\right)\Lambda^2(n)\nonumber\\
&=\sum_{L\in \mathcal{L}_1}Q_1(L)-Q_2-\sum_{L\in\mathcal{L}_2}Q_3(L),\label{eq:SSplit}
\end{align}
where
{\allowdisplaybreaks
\begin{align}
Q_1(L)&=\sum_{N\le n\le 2N}\chi_1(L(n))\Lambda^2(n),\\
Q_2&=\sum_{N\le n\le 2N}\Lambda^2(n),\\
Q_3(L)&=\sum_{N\le n\le 2N}\sum_{\substack{p|L(n)\\p\le (2N)^{1/(r+1)}}}\Lambda^2(n).
\end{align}}
The choice of good values for $\lambda_d$ and the corresponding evaluation of $Q_1$, $Q_2$, $Q_3$ already exists in the literature. We quote from \cite{Maynard}[Proposition 4.1] taking $W_0(t)=1$ and\cite{GGPY}[Theorem 7 and Theorem 9]. We note that \cite{GGPY}[Theorem 9] does not require that $E_2$-numbers have level of distribution $\theta$, so Hypothesis $BV(\theta)$ is sufficient for the statement to hold. These results give for a fixed polynomial $P$, for $k=k_1+k_2=\#(\mathcal{L}_1\cup\mathcal{L}_2)$, for $L\in\mathcal{L}$ and for sufficiently large $C$ that we can choose the $\lambda_d$ such that
\begin{align}
Q_1(L)&\sim\frac{\mathfrak{S}(\mathcal{L})N(\log{R})^{k+1}}{(k-2)!\log{N}}\int_0^1\tilde{P}(1-t)^2t^{k-2}d t,\\
Q_2&\sim\frac{\mathfrak{S}(\mathcal{L})N(\log{R})^{k}}{(k-1)!}\int_0^1P(1-t)^2t^{k-1}d t,\\
Q_3(L)&\sim\frac{\mathfrak{S}(\mathcal{L})N(\log{R})^{k}}{(k-1)!}I,\label{eq:Q3Def}
\end{align}
where
\begin{align}
\mathfrak{S}(\mathcal{L})&\text{ is a positive constant depending only on $\mathcal{L}$,}\\
\tilde{P}(x)&=\int_0^xP(t)d t,\\
I&=\int_0^{\delta}\frac{1}{y}\int_{1-y}^1P(1-t)^2t^{k-1}d t d y\nonumber\\
&\qquad+\int_0^{\delta}\frac{1}{y}\int_0^{1-y}\left(P(1-t)-P(1-t-y)\right)^2t^{k-1}d td y,\\
\delta&=\frac{2}{\theta(r+1)}.
\end{align}
Here the asymptotic for $Q_3$ is valid only for $R^2(2N)^{1/(r+1)}\le N(\log{N})^{-C}$, and so we introduce the condition
\begin{equation}\label{eq:rCondition}
r+1>\frac{1}{1-\theta}
\end{equation}
to ensure that this is satisfied for $N$ sufficiently large. All the other asymptotics are valid without further conditions.

We choose the polynomial 
\begin{equation}
P(x)=x^l,
\end{equation}
 where $l\ge 0$ is an integer to be declared later. To ease notation we let
\begin{equation}
C(\mathcal{L})=\frac{\mathfrak{S}(\mathcal{L})N(\log{R})^k(2l)!}{(k+2l)!}.\label{eq:CDef}
\end{equation}
We note that for $a,b\in \mathbb{N}$ 
\begin{equation}
\int_0^1x^a(1-x)^b d x =\frac{a!b!}{(a+b+1)!}.
\end{equation}
Thus we see that
\begin{align}
\sum_{L\in\mathcal{L}_1}Q_1(L)&\sim \theta\left(\frac{2l+1}{l+1}\right)\left(\frac{k_1}{k+2l+1}\right)C(\mathcal{L}),\label{eq:Q1}\\
Q_2&\sim C(\mathcal{L}).\label{eq:Q2}
\end{align}
We follow a similar approach to Graham, Goldston, Pintz and Y\i ld\i r\i m  \cite{GGPY} to estimate $Q_3$. We let
\begin{equation}
I=\int_0^\delta \frac{F(y)}{y}d y,\label{eq:Idef}
\end{equation}
where
\begin{align}
F(y)&=\int_{1-y}^1 P(1-t)^2t^{k-1}d t+\int_0^{1-y}\left(P(1-t)-P(1-t-y)\right)^2t^{k-1}d t\nonumber\\
&=\int_0^1P(1-t)^2t^{k-1}d t+\int_0^{1-y}P(1-t-y)^2t^{k-1}d t\nonumber\\
&\qquad-2\int_0^{1-y}P(1-t)P(1-t-y)t^{k-1}d t.
\end{align}
We recall that $P(x)=x^l$, and note that
\begin{equation}
P(1-t)=(1-t)^l=(1-t-y)^l+\sum_{j=1}^l\binom{l}{j}y^j(1-t-y)^{l-j}.
\end{equation}
Thus
\begin{align}
F(y)&=\int_0^1(1-t)^{2l}t^{k-1}d t+\int_0^{1-y}(1-t-y)^{2l}t^{k-1}d t\nonumber\\
&\qquad-2\int_0^{1-y}(1-t-y)^{2l}t^{k-1}d t-2\sum_{j=1}^l\binom{l}{j}y^j\int_0^{1-y}(1-t-y)^{2l-j}t^{k-1}d t\nonumber\\
&\le \int_0^1(1-t)^{2l}t^{k-1}d t-\int_0^{1-y}(1-t-y)^{2l}t^{k-1}d t\nonumber\\
&=\frac{(k-1)!(2l)!}{(k+2l)!}\left(1-(1-y)^{k+2l}\right).
\end{align}
Substituting this into \eqref{eq:Idef} gives
\begin{align}
I=\int_0^\delta\frac{F(y)}{y}d y &\le \frac{(k-1)!(2l)!}{(k+2l)!}\int_0^\delta\frac{1-(1-y)^{k+2l}}{y}d y\\
&=\frac{(k-1)!(2l)!}{(k+2l)!}\sum_{j=0}^{k+2l-1}\int_0^\delta (1-y)^j dy\\
&\le \frac{(k-1)!(2l)!}{(k+2l)!}\delta(k+2l).\label{eq:IVal}
\end{align}
Therefore from \eqref{eq:Q3Def}, \eqref{eq:CDef} and \eqref{eq:IVal}, for any $L\in\mathcal{L}_2$ we have that
\begin{equation}\label{eq:Q3}
Q_3(L)\le C(\mathcal{L})\left(\delta(k+2l)+o(1)\right).
\end{equation}
Substituting \eqref{eq:Q1}, \eqref{eq:Q2} and \eqref{eq:Q3} into \eqref{eq:SSplit} we see
\begin{equation}
S\ge C(\mathcal{L})\left(\theta\left(\frac{2l+1}{l+1}\right)\left(\frac{k_1}{k+2l+1}\right)-1-k_2\delta(k+2l)+o(1)\right).
\end{equation}
We let
\begin{equation}
k+2l+1=\lceil C_1(2\theta-1)^{-2}\rceil ,\qquad l+1=\lceil C_2(2\theta-1)^{-1}\rceil\label{eq:klDef}
\end{equation}
for some $C_1$, $C_2$. This gives
\begin{align}
\frac{S}{C(\mathcal{L})}&\ge \theta\left(2-\frac{2\theta-1}{C_2}\right)\left(1-\frac{2C_2(2\theta-1)}{C_1}-\frac{(k_2+1)(2\theta-1)^2}{C_1}\right)-1\nonumber\\
&\qquad-k_2\delta C_1(2\theta-1)^{-2}+o(1)\nonumber\\
&=(2\theta-1)\left(1-\frac{\theta}{C_2}-\frac{4\theta C_2}{C_1}-\frac{2k_2(2\theta-1)\theta}{C_1}+\frac{(1+k_2)\theta(2\theta-1)^2}{C_1C_2}\right)\nonumber\\
&\qquad-k_2\delta C_1(2\theta-1)^{-2}+o(1).
\end{align}
We let
\begin{equation}
C_1=40k_2,\qquad C_2=3.
\end{equation}
We see from \eqref{eq:klDef} that this choice of $C_1$ and $C_2$ corresponds to positive integer values for $k_1$ and $l$ for any choice of $0.5<\theta\le 0.99$ or $k_2$, and so is a valid choice.

Since $k_2$ is a positive integer and $1/2<\theta\le 1$, this gives
\begin{align}
\frac{S}{C(\mathcal{L})}&\ge (2\theta-1)\left(\frac{30-17\theta-5\theta^2+2\theta^3}{30}\right)-40k_2^2(2\theta-1)^{-2}\delta+o(1)\nonumber\\
&\ge (2\theta-1)\left(\frac{31-21\theta}{30}\right)-40k_2^2(2\theta-1)^{-2}\delta+o(1).
\end{align}
Thus $S>0$ for large $N$ if $\delta$ is chosen such that
\begin{equation}
\delta < \frac{(2\theta-1)^3(31-21\theta)}{1200k_2^2}.
\end{equation}
We recall $\delta=2/\theta(r+1)$, so $S$ is positive provided $r$ is chosen larger than
\begin{equation}
\frac{2400k_2^2}{(2\theta-1)^3\theta(31-21\theta)}-1< \frac{240k_2^2}{(2\theta-1)^3}.
\end{equation}
We note that if $r=240k_2^2/(2\theta-1)^3$ then for $\theta<0.99$ the condition \eqref{eq:rCondition} is satisfied. This completes the proof of Theorem \ref{thrm:MainTheorem}.

We remark that by choosing $L(n)=n+h$ with $h<H$ for $L\in\mathcal{L}_2$ and $L(n)=n+h$ with $h>H$ for $L\in\mathcal{L}_1$ we can ensure that of the $k_2+2$ almost-primes we find, the largest two are primes.
\section{Proof of Theorem \ref{thrm:ElliottHalberstam}}
We can get better quantitative results for the number of prime factors involved in our almost-prime if we assume a fixed level of distribution result for almost-primes and for primes, and then follow the work of \cite{GGPY}.

We consider the same sum $S$, but now we assume that $\mathcal{L}_2=\{L_0\}$ and we take $r=4$. Thus $k_2=1$ and $k=k_1+1$.
\begin{align}
S=S(N;\mathcal{L}_1,\{h_0\})&=\sum_{N\le n \le 2N}\left(\sum_{L\in\mathcal{L}_1}\chi_{1}(L(n))+\chi_{4}(L_0(n))-2\right)\left(\sum_{\substack{d|\Pi(n)\\ d\le R}}\lambda_d\right)^2\nonumber\\
&=\sum_{L\in\mathcal{L}_1}Q_1(L)+Q_1'(L_0)-Q_2,
\end{align}
where $Q_1(L)$, $Q_2$ are as before and
\begin{align}
Q_1'&=\sum_{N\le n \le 2N}\chi_{4}(L_0(n))\left(\sum_{\substack{d|\Pi(n)\\ d\le R}}\lambda_d\right)^2\\
R&=N^{0.99/2}(\log{N})^C.
\end{align}
As before, $C$ is a suitably large positive constant.

We split the contribution to $Q_1'$ depending on whether $L_0(n)$ has exactly 1, 2, 3 or 4 prime factors. Thus
\begin{equation}
Q_1'=Q_{1}(L_0)+Q_{12}'+Q_{13}'+Q_{14}',
\end{equation}
where
\begin{equation}
Q_{1j}'=\sum_{N\le n\le 2N}\beta_j(L_0(n))\left(\sum_{\substack{d|\Pi(n)\\ d\le R}}\lambda_d\right)^2,
\end{equation}
and
\begin{equation}
\beta_j(n)=\begin{cases}
1,\qquad &\text{$n$ has exactly $j$ prime factors}\\
0, &\text{otherwise.}
\end{cases}
\end{equation}
For technical reasons we find it harder to deal with terms arising when $L(n)$ has a prime factor less than $N^{\epsilon}$ or no prime factor greater than $N^{1/2}$. Thus we obtain a lower bound for $Q_{1j}'$ by replacing $\beta_j(L_0(n))$ with $\beta'_j(L_0(n))$, where
\begin{equation}
\beta_j'(n)=\begin{cases}
1,\qquad &\text{$n=p_1p_2\dots p_j$ with $n^\epsilon<p_1<\dots<p_j$ and $n^{0.505}<p_j$}\\
0, &\text{otherwise.}
\end{cases}
\end{equation}
We can then obtain these asymptotic lower bounds. By following an equivalent argument to \cite{Thorne} and \cite{Maynard}[Proposition 4.2] but using Hypothesis GBV(0.99,$j$) to bound the error terms we have
\begin{equation}
Q_{1j}'\ge (1+o(1))\frac{\mathfrak{S}(\mathcal{L})(\log{R})^{k+1}}{(k-2)!\log{N}}J_{j},
\end{equation}
where
{\allowdisplaybreaks
\begin{align}
J_{r}&=\int_{(x_1,\dots,x_{r-1})\in \mathcal{A}_r}\frac{ I_1(Bx_1,\dots,Bx_{r-1})}{\left(\prod_{i=1}^{r-1}x_i\right)\left(1-\sum_{i=1}^{r-1}x_i\right)}d x_1\dots dx_{r-1},\\
I_1&=\int_0^1\left(\sum_{J\subset\{1,\dots,r-1\}}(-1)^{|J|}\tilde{P}^+(1-t-\sum_{i\in J}x_i) \right)^2t^{k-2}d t,\\
\tilde{P}^+(x)&=\begin{cases}
\int_0^x P(t)d t,\qquad &x\ge 0\\
0,&\text{otherwise,}
\end{cases}\\
B&=\frac{2}{0.99},\\
\mathcal{A}_r&=\left\{x\in[0,1]^{r-1}:\epsilon<x_1<\dots<x_{r-1},\sum_{i=1}^{r-1}x_i<B^{-1}\right\}.
\end{align}}
As before, by Hypothesis BV(0.99) we also have for any $L\in\mathcal{L}$
\begin{align}
Q_1(L)&\sim\frac{\mathfrak{S}(\mathcal{L})(\log{R})^{k+1}}{(k-2)!\log{N}}\int_0^1\tilde{P}(1-t)^2t^{k-2}d t,\\
Q_2&\sim\frac{\mathfrak{S}(\mathcal{L})(\log{R})^{k}}{(k-1)!}\int_0^1P(1-t)^2t^{k-1}d t.
\end{align}
Thus we have that
\begin{equation}
S\ge \frac{\mathfrak{S}(\mathcal{L})(\log{R})^{k}}{(k-1)!}\left(\frac{0.99(k-1)}{2}\left(k J_1+J_2+J_3+J_4\right)-2I_0+o(1)\right),
\end{equation}
where $J_r$ is given above and
\begin{equation}
I_0=\int_0^1P(1-t)^2t^{k-1}d t.
\end{equation}
Therefore given a polynomial $P$ we can get an asymptotic lower bound for $S$ by explicitly calculating the integrals $I_0$, $J_1$, $J_2$, $J_3$ and $J_4$.

Explicitly we have for $r=1$
\begin{equation}
J_1=\int_0^1\tilde{P}(1-t)^2t^{k-2}d t.
\end{equation}
Similarly for $r=2$ we have
\begin{equation}
J_2=J_{21}+J_{22}+O(\epsilon),
\end{equation}
where
\begin{align}
J_{21}&=\int_{0}^{1}\frac{B}{x(B-x)}\int_0^{1-x}\left(\tilde{P}(1-t)-\tilde{P}(1-t-x)\right)^2t^{k-2}d t d x,\\
J_{22}&=\int_{0}^{1}\frac{B}{x(B-x)}\int_{1-x}^1\tilde{P}(1-t)^2t^{k-2}d t d x.
\end{align}
Similarly for $r=3$ we have
\begin{equation}
J_3=J_{31}+J_{32}+J_{33}+J_{34}+O(\epsilon),
\end{equation}
where
{\allowdisplaybreaks
\begin{align}
J_{31}&=\int_{0}^{1/2}\int_{x}^{1-x}\frac{B}{xy(B-x-y)}\int_{1-x}^1\left(\tilde{P}(1-t)\right)^2t^{k-2}d t d y d x,\\
J_{32}&=\int_{0}^{1/2}\int_{x}^{1-x}\frac{B}{xy(B-x-y)}\int_{1-y}^{1-x}\left(\tilde{P}(1-t)-\tilde{P}(1-t-x)\right)^2t^{k-2}d t d y d x,\\
J_{33}&=\int_{0}^{1/2}\int_{x}^{1-x}\frac{B}{xy(B-x-y)}\int_{1-x-y}^{1-y}\nonumber\\
&\quad\left(\tilde{P}(1-t)-\tilde{P}(1-t-x)-\tilde{P}(1-t-y)\right)^2t^{k-2}d t d y d x,\\
J_{34}&=\int_{0}^{1/2}\int_{x}^{1-x}\frac{B}{xy(B-x-y)}\int_{0}^{1-x-y}\nonumber\\
&\quad\left(\tilde{P}(1-t)-\tilde{P}(1-t-x)-\tilde{P}(1-t-y)+\tilde{P}(1-t-x-y)\right)^2t^{k-2}d t d y d x.
\end{align}}
Finally for $r=4$ we have
\begin{equation}
J_4=J_{41}+J_{42}+J_{43}+J_{44}+J_{45}+J_{46}+J_{47}+J_{48}+J_{49}+J_{410}+J_{411}+O(\epsilon),
\end{equation}
where
{\allowdisplaybreaks
\begin{align}
J_{41}&=\int_{0}^{1/3}\int_{x}^{(1-x)/2}\int_{y}^{1-x-y}\frac{B}{xyz(B-x-y-z)}\int_{1-x}^1\left(\tilde{P}(1-t)\right)^2t^{k-2}d t d z d y d z,\\
J_{42}&=\int_{0}^{1/3}\int_{x}^{(1-x)/2}\int_{y}^{1-x-y}\frac{B}{xyz(B-x-y-z)}\int_{1-y}^{1-x}\nonumber\\
&\qquad\left(\tilde{P}(1-t)-\tilde{P}(1-t-x)\right)^2t^{k-2}d t d z d y d z,\\
J_{43}&=\int_{0}^{1/4}\int_{x}^{1/2-x}\int_{x+y}^{1-x-y}\frac{B}{xyz(B-x-y-z)}\int_{1-x-y}^{1-y}\nonumber\\
&\qquad\left(\tilde{P}(1-t)-\tilde{P}(1-t-x)-\tilde{P}(1-t-y)\right)^2t^{k-2}d t d z d y d z,\\
J_{44}&=\int_{0}^{1/3}\int_{x}^{(1-x)/2}\int_{y}^{x+y}\frac{B}{xyz(B-x-y-z)}\int_{1-z}^{1-y}\nonumber\\
&\qquad\left(\tilde{P}(1-t)-\tilde{P}(1-t-x)-\tilde{P}(1-t-y)\right)^2t^{k-2}d t d z d y d z,\\
J_{45}&=\int_{0}^{1/4}\int_{x}^{1/2-x}\int_{x+y}^{1-x-y}\frac{B}{xyz(B-x-y-z)}\int_{1-z}^{1-x-y}\nonumber\\
&\quad\left(\tilde{P}(1-t)-\tilde{P}(1-t-x)-\tilde{P}(1-t-y)+\tilde{P}(1-t-x-y)\right)^2t^{k-2}d t d z d y d z,\\
J_{46}&=\int_{0}^{1/3}\int_{x}^{(1-x)/2}\int_{y}^{x+y}\frac{B}{xyz(B-x-y-z)}\int_{1-x-y}^{1-z}\nonumber\\
&\quad\left(\tilde{P}(1-t)-\tilde{P}(1-t-x)-\tilde{P}(1-t-y)-\tilde{P}(1-t-z)\right)^2t^{k-2}d t d z d y d z,\\
J_{47}&=\int_{0}^{1/4}\int_{x}^{1/2-x}\int_{x+y}^{1-x-y}\frac{B}{xyz(B-x-y-z)}\int_{1-x-z}^{1-z}\nonumber\\
&\quad\Bigl(\tilde{P}(1-t)-\tilde{P}(1-t-x)-\tilde{P}(1-t-y)\nonumber\\
&\quad\quad-\tilde{P}(1-t-z)+\tilde{P}(1-t-x-y)\Bigr)^2t^{k-2}d t d z d y d z,\\
J_{48}&=\int_{0}^{1/3}\int_{x}^{(1-x)/2}\int_{y}^{x+y}\frac{B}{xyz(B-x-y-z)}\int_{1-x-z}^{1-x-y}\nonumber\\
&\quad\Bigl(\tilde{P}(1-t)-\tilde{P}(1-t-x)-\tilde{P}(1-t-y)-\tilde{P}(1-t-z)\nonumber\\
&\quad\quad+\tilde{P}(1-t-x-y)\Bigr)^2t^{k-2}d t d z d y d z,\\
J_{49}&=\int_{0}^{1/3}\int_{x}^{(1-x)/2}\int_{y}^{1-x-y}\frac{B}{xyz(B-x-y-z)}\int_{1-y-z}^{1-x-z}\nonumber\\
&\quad\Bigl(\tilde{P}(1-t)-\tilde{P}(1-t-x)-\tilde{P}(1-t-y)-\tilde{P}(1-t-z)\nonumber\\
&\quad\quad+\tilde{P}(1-t-x-y)+\tilde{P}(1-t-x-z)\Bigr)^2t^{k-2}d t d z d y d z,\\
J_{410}&=\int_{0}^{1/3}\int_{x}^{(1-x)/2}\int_{y}^{1-x-y}\frac{B}{xyz(B-x-y-z)}\int_{1-x-y-z}^{1-y-z}\nonumber\\
&\quad\Bigl(\tilde{P}(1-t)-\tilde{P}(1-t-x)-\tilde{P}(1-t-y)-\tilde{P}(1-t-z)+\tilde{P}(1-t-x-y)\nonumber\\
&\quad\quad+\tilde{P}(1-t-x-z)+\tilde{P}(1-t-y-z)\Bigr)^2t^{k-2}d t d z d y d z,\\
J_{411}&=\int_{0}^{1/3}\int_{x}^{(1-x)/2}\int_{y}^{1-x-y}\frac{B}{xyz(B-x-y-z)}\int_0^{1-x-y-z}\nonumber\\
&\Bigl(\tilde{P}(1-t)-\tilde{P}(1-t-x)-\tilde{P}(1-t-y)-\tilde{P}(1-t-z)+\tilde{P}(1-t-x-y)\nonumber\\
&+\tilde{P}(1-t-x-z)+\tilde{P}(1-t-y-z)-\tilde{P}(1-t-x-y-z)\Bigr)^2t^{k-2}d t d z d y d z.
\end{align}}
We choose $k=22$ and $P(t)=1+60t-300t^2+3500t^3$ and find that
\begin{align}
I_0&=\frac{121351}{59202}=2.04978\dots,\\
J_1&=\frac{228380}{18027009}=0.01266\dots,\\
J_2&\ge 0.041+O(\epsilon),\\
J_3&\ge 0.048+O(\epsilon),\\
J_4&\ge 0.028+O(\epsilon).
\end{align}
Thus we have that
\begin{align}
S&\ge \frac{\mathfrak{S}(\mathcal{L})(\log{R})^{k}}{(k-1)!}\left(\frac{0.99(k-1)}{2}\left(k J_1+J_2+J_3+J_4\right)-2I_0+O(\epsilon)+o(1)\right)\nonumber\\
&\ge \frac{\mathfrak{S}(\mathcal{L})(\log{R})^{k}}{(k-1)!}\left(0.013+O(\epsilon)+o(1)\right).
\end{align}
In particular, for $N$ sufficiently large and $\epsilon$ sufficiently small we have $S>0$, and so there are infinitely many $n$ for which an admissible 22-tuple attains at least two prime values and one value with at most 4 prime factors.

The set $\{0, 6, 8, 14, 18, 20, 24, 30, 36, 38, 44, 48, 50, 56, 60, 66, 74, 78, 80, 84, 86, 90\}$ is an admissible 22-tuple, and so the interval $[n,n+90]$ infinitely often contains at least two primes and an integer with at most 4 prime factors.

We remark here that if we can take the level of distribution $\theta=1-\delta$ for every $\delta>0$ then we can take $k=19$ instead of 22, which reduces the length of the interval to 80.

\section{Acknowledgment}
I would like to thank my supervisor, Prof. Heath-Brown, for suggesting this problem and for his careful reading of this paper.
\bibliographystyle{acm}
\bibliography{bibliography}
\end{document}